\documentclass[11pt,a4paper]{article}

\usepackage{amsmath,amssymb,}

\newtheorem{thm}{Theorem}[section]
\newtheorem{lem}[thm]{Lemma}
\newtheorem{prop}[thm]{Proposition}

\newtheorem{cor}[thm]{Corollary}
\newtheorem{ex}[thm]{Examples}
\newtheorem{rem}[thm]{Remark}

\begin{document}

\title{\huge \textbf{On central sequence algebras of tensor product von Neumann algebras}}
\author{Yasuhito Hashiba}
\date{}
\maketitle

\begin{abstract}
\noindent We show that when $M,N_{1},N_{2}$ are tracial von Neumann algebras with $M'\cap M^{\omega}$ abelian, $M'\cap(M\bar{\otimes}N_{1})^{\omega}$ and  $M'\cap(M\bar{\otimes}N_{2})^{\omega}$ commute in  $(M\bar{\otimes}N_{1}\bar{\otimes}N_{2})^{\omega}$. As a consequence, we obtain information on McDuff decomposition of $\rm{I\hspace{-.01em}I}_{1}$ factors  of the form $M\bar{\otimes}N$, where $M$ is a non-McDuff factor.
\end{abstract}

\section{INTRODUCTION} 
When we have a tracial von Neumann algebra $M$, the \textit{central sequence algebra} is defined as $M'\cap M^{\omega}$, where  
$\omega$ is a free ultrafilter on $\mathbb{N}$ and $M^{\omega}$ is the ultraproduct von Neumann algebra \cite{mcduff1970central}. Central sequences and  central sequence algebras (or more generally, $N'\cap M^{\omega}$ when we have $N \subset M$) have been very important in the study of von Neumann algebras. In \cite{murray1943rings}, Murray and von Neumann introduced \textit{property Gamma} for  $\rm{I\hspace{-.01em}I}_{1}$ factors (which is equivalent to  $M'\cap M^{\omega}\neq \mathbb{C}$ in terms of ultraproducts). They showed that the hyperfinite  $\rm{I\hspace{-.01em}I}_{1}$ factor $R$ has property Gamma, while the free group factor $L(\mathbb{F}_{2})$ does not, thus giving the first example of two non-isomorphic separable  $\rm{I\hspace{-.01em}I}_{1}$ factors. Later on, McDuff showed that $M\cong M\bar{\otimes}R$ and $M'\cap M^{\omega}$ being non-abelian are equivalent \cite{mcduff1970central}. Factors of type $\rm{I\hspace{-.01em}I}_{1}$ satisfying these conditions are now called \textit{McDuff factors}.\vspace{\baselineskip}

The tensor product construction is another natural notion in von Neumann algebra theory which have been studied over the years. Therefore, it is natural to consider ultraproducts and central sequence algebras of von Neumann algebras which is a tensor product of two or more von Neumann algebras. In this direction, we have that $R'\cap R^{\omega}$ is a prime  $\rm{I\hspace{-.01em}I}_{1}$ factor (i.e. it can not be decomposed into a tensor product of two  $\rm{I\hspace{-.01em}I}_{1}$ factors), while we have $R\cong R\bar{\otimes}R$. This observation shows that $(M\bar{\otimes}N)'\cap (M\bar{\otimes}N)^{\omega}$ can be strictly larger than $(M'\cap M^{\omega})\bar{\otimes} (N'\cap N^{\omega})$. Note that when $M$ is a \textit{full factor}, i.e. $M'\cap M^{\omega}=\mathbb{C}$,  we have $(M\bar{\otimes}N)'\cap (M\bar{\otimes}N)^{\omega}= N'\cap N^{\omega}$. See \cite{fang2006central} for the facts in this paragraph.\vspace{\baselineskip}

In this paper, we continue the investigation of the relationship between central sequence algebras and tensor products. Our main result is the following theorem.\vspace{\baselineskip}

\textbf{Theorem \ref{theorem}} \textit{ Let $M,N_{1},N_{2}$ be tracial von Neumann algebras. Assume that $M'\cap M^{\omega}$ is abelian.
Then $M'\cap(M\bar{\otimes}N_{1})^{\omega}$ and  $M'\cap(M\bar{\otimes}N_{2})^{\omega}$ commute in  $(M\bar{\otimes}N_{1}\bar{\otimes}N_{2})^{\omega}$. }\vspace{\baselineskip}

This theorem applies when $M$ is a  non McDuff factor. If we can additionally assume that  $M'\cap(M\bar{\otimes}N_{1})^{\omega}=(M'\cap M^{\omega})\bar{\otimes}N_{1}^{\omega}$ and $M'\cap(M\bar{\otimes}N_{2})^{\omega}=(M'\cap M^{\omega})\bar{\otimes}N_{2}^{\omega} $ (this is always satisfied when $M$ is a full factor), our theorem follows directly from the definitions. However, similar to the fact explained in the second paragraph, $M'\cap (M\bar{\otimes}N)^{\omega}$ can be strictly larger than $(M'\cap M^{\omega})\bar{\otimes}N^{\omega}$ (see Lemma \ref{lem2}). Our theorem still holds in such a general case.\vspace{\baselineskip}

The proof of the theorem uses ideas of \cite{marrakchi2018stability} and  \cite{ioana2015spectral}. (See also \cite{haagerup1985new}.) A key step of the proof is a reduction of the problem to one in a much easier case: when $N_{1}$ and $N_{2}$ is abelian. The easier case is proved in Proposition \ref{proposition}. As a corollary, we obtain information on 
McDuff decomposition of $\rm{I\hspace{-.01em}I}_{1}$ factors  of the form $M\bar{\otimes}N$, where $M$ is a non-McDuff factor (see Corollary \ref{cor3.4} and Corollary \ref{cor3.5}).

\section{PRELIMINARIES}  
\subsection{Terminology}
A tracial von Neumann algebra $(M, \tau)$ is a von Neumann algebras $M$ endowed with a faithful normal tracial state $\tau : M \to \mathbb{C}$. We denote by $L^{2}(M)$ the Hilbert space obtained by completion of $M$ with respect to the norm $\parallel x \parallel_{2}= \tau (x^{*}x)^{1/2}$. We say that $M$ is separable if it is separable with respect to this norm. Every von Neumann algebra considered in this paper is assumed to be separable (except for ultraproducts). We denote by $\mathcal{U}(M)$ the group of unitaries of $M$. For $r\ge 0$, we set $(M)_{r}= \{x\in M | \parallel x \parallel \le r \}$. For a von Neumann subalgebra $N \subset M$, we regard $N$ as a tracial von Neumann algebra together with $\tau |_{N}$, and we denote by $E_{N}$ the trace preserving normal conditional expectation $M \to N$. For von Neumann subalgebras $K,L \subset M$, we denote by $K \vee L$ the von Neumann algebra generated by $K$ and $L$. The hyperfinite  $\rm{I\hspace{-.01em}I}_{1}$ factor is denoted as $R$. We have $R \cong \bar{\otimes}_{n}^{\infty} M_{2}(\mathbb{C})$. We refer to \cite{anantharaman2017introduction} for the fundamental facts about (tracial) von Neumann algebras used in this paper.

\subsection{Ultraproducts}
Throughout this paper, we fix a free ultrafilter $\omega$ on $\mathbb{N}$. 
For any sequence $(M_{n}, \tau_{n})$ of tracial von Neumann algebras, we define the ultraproduct algebra $\Pi_{\omega} M_{n}$ as the quotient  $\textrm C^{*}$-algebra $A/I$, where $A=\{(x_{n})_{n} \in \Pi_{n} M_{n} | \sup_{n} \parallel x_{n} \parallel < \infty \}$ and $I=\{ (x_{n})_{n} \in A | \lim_{n \to \omega} \parallel x_{n} \parallel_{2} =0\}$. This is a tracial von Neumann algebra with trace $\tau ((x_{n})_{n})= \lim_{n \to \omega} \tau (x_{n})$. When we have $M_{n}=M$ for all $n$, we denote the ultraproduct algebra by $M^{\omega}$. In this case, $M$ is naturally embedded into $M^{\omega}$ by $x \in M \mapsto (x)_{n} \in M^{\omega}$. We regard $M$ as a tracial von Neumann subalgebra of $M^{\omega}$ by this embedding.

\subsection{Central sequence algebras}
For a tracial von Neumann algebra $M$, the relative commutant $M'\cap M^{\omega}$ is called the central sequence algebra of $M$. This is a very useful tool to distinguish  $\rm{I\hspace{-.01em}I}_{1}$ factors.

\begin{thm} (\cite{mcduff1970central,dixmiercentral})
Let $M$ be a (separable) $\rm{I\hspace{-.01em}I}_{1}$ factor. Then $M$ satisfies one (and only one) of the following conditions.
\begin{enumerate}
\item
We have $M' \cap M^{\omega}= \mathbb{C}$. In this case, we say that $M$ is a full factor.
\item
We have that $M' \cap M^{\omega}$ is diffuse and abelian.
\item
We have that $M' \cap M^{\omega}$ is a type  $\rm{I\hspace{-.01em}I}_{1}$ von Neumann algebra. This is equivalent to the following conditions.
 \begin{enumerate}
  \item
  $M' \cap M^{\omega}$ is not abelian.
  \item
  $M \cong M \bar{\otimes} R$, where $R$ is the hyperfinite $\rm{I\hspace{-.01em}I}_{1}$ factor.
  \end{enumerate}
  In this case, we say that $M$ is a McDuff factor.
  We say that a McDuff factor $M$ admits a McDuff decomposition if it can be written as $M=N\bar{\otimes}R$ for some non-McDuff $\rm{I\hspace{-.01em}I}_{1}$ factor $N$.
\end{enumerate}
For $M$, satisfying property Gamma is equivalent to not being a full factor in the above sense.
\end{thm}

\begin{ex}
\begin{itemize}
\item
The free group factor $L(\mathbb{F}_{n})$ is a full factor for all $n \ge 2$. More generally, by \cite{effros1975property}, if $G$ is a discrete i.c.c. group which is not inner amenable, then $L(G)$ is a full factor. (The converse does not hold: there is an i.c.c. group which is inner amenable with $L(G)$ a full factor(\cite{vaes2012inner}).)
\item
The hyperfinite  $\rm{I\hspace{-.01em}I}_{1}$ factor is McDuff. More generally, if  $(M_{n})_{n}$ is a sequence of nontrivial tracial factors, then the infinite tensor product $\bar{\otimes}_{n}^{\infty} M_{n}$ is a McDuff  $\rm{I\hspace{-.01em}I}_{1}$ factor.
\end{itemize}
\end{ex}

We end this section with an easy lemma. It is probably well known for experts, but we include a proof for completeness. It is not used in the proof of the main theorem, but it shows that the statement of the theorem is not trivial. 
\begin{lem}
Let $M, N$ be tracial von Neumann algebras. Assume that $M'\cap M^{\omega}$ and $N$ are diffuse. Then $M'\cap (M\bar{\otimes}N)^{\omega} $ is strictly larger than $(M'\cap M^{\omega})\bar{\otimes}N^{\omega}$.  \label{lem2}
\end{lem}

\textit{Proof.} It is easy to see that there is a natural inclusion $(M'\cap M^{\omega})\bar{\otimes}N^{\omega}\subset M'\cap (M\bar{\otimes}N)^{\omega}$ , so we only show that they do not coincide.\vspace{\baselineskip}

Take $(F_{n})_{n}$ as an increasing sequence of finite sets of $M$ such that $\cup_{n} F_{n}$ is dense in $M$. Since $M'\cap M^{\omega}$ is diffuse, we can take a Haar unitary $u\in M'\cap M^{\omega}$. Writing $u^{i} $ by sequences in $M$ for  $ i\in \mathbb{N}$, for each $n\in \mathbb{N}$, we can take an orthogonal family $\{y_{i}^{n}\}_{i}\in (M)_{1}$ such that we have $\sum_{v\in F_{n}}\parallel [v,y_{i}^{n}] \parallel_{2}^{2} <\frac{1}{n} $ and $\parallel y_{i}^{n} \parallel_{2}^{2}\ge \frac{1}{2}$ for each $ 1\le i \le n $.  \vspace{\baselineskip}

Since $N$ is diffuse, we can take $A\subset N$ with $A\cong L^{\infty}([0,1])$. For each $n$, set $x_{n}=\sum_{i=1}^{n} y_{i}^{n} \otimes \chi_{[\frac{i-1}{n},\frac{i}{n})} \in (M\bar{\otimes}A)_{1}$, where $\chi_{[\frac{i-1}{n},\frac{i}{n})} $ is the characteristic function in $A\cong L^{\infty}([0,1])$.
For each $n$, we have 
\begin{align*}
\sum_{v\in F_{n}} \parallel [v,x_{n}] \parallel_{2}^{2} 
 &\le \sum_{F_{n}} \int_{[0,1]} \sum_{i=1}^{n} \parallel [v,y_{i}^{n}]\parallel_{2}^{2} \chi_{[\frac{i-1}{n},\frac{i}{n})}(t) dt \\
 &\le \frac{1}{n} \sum_{i}\int   \chi_{[\frac{i-1}{n},\frac{i}{n})}(t) dt \\
 &= \frac{1}{n}.
\end{align*}
This shows that $x=(x_{n})_{n}\in M'\cap (M\bar{\otimes}A)^{\omega}$. Note that $\parallel x \parallel_{2}^{2}\ge\frac{1}{2}$. \vspace{\baselineskip}

We show that this $x$ is orthogonal to $M^{\omega}\bar{\otimes}A^{\omega}$. Take any $(z_{n})_{n}\in (M^{\omega})_{1}$ and $(a_{n})_{n}\in (A^{\omega})_{1}$ with $\parallel z_{n} \parallel, \parallel a_{n} \parallel \le 1$ for all $n\in \mathbb{N}$. We have 
\begin{align*}
|\langle z_{n}\otimes a_{n}, x_{n}\rangle| &= |\int \sum_{i=1}^{n} \langle z_{n}, y_{i}^{n}\rangle \overline{a_{n}(t)}\chi_{[\frac{i}{n},\frac{i+1}{n})}(t) dt |\\
&\le \sum_{i} |\langle z_{n},y_{i}^{n}\rangle | \int|\overline{a_{n}(t)}\chi_{[\frac{i}{n},\frac{i+1}{n})}(t)| dt\\
&\le (\sum_{i}  |\langle z_{n},y_{i}^{n}\rangle |^{2})^{\frac{1}{2}} (\sum_{i} (\int|\overline{a_{n}(t)}\chi_{[\frac{i}{n},\frac{i+1}{n})}(t)| dt)^{2})^{\frac{1}{2}}.
\end{align*}
Since  $\{y_{i}^{n}\}_{i}$ is orthogonal, we have $\sum_{i}  |\langle z_{n},y_{i}^{n}\rangle |^{2}\le \parallel z_{n} \parallel_{2}^{2}\le 1$. Also, since $|a_{n}(t)|\le 1$, we can compute 
$\sum_{i} (\int|\overline{a_{n}(t)}\chi_{[\frac{i}{n},\frac{i+1}{n})}(t)| dt)^{2}\le \sum_{i} (\int \chi_{[\frac{i}{n},\frac{i+1}{n})}(t) dt)^{2}\ \le \sum_{i=1}^{n} (\frac{1}{n})^{2} \le \frac{1}{n}$. Hence, we have $|\langle z_{n}\otimes a_{n}, x_{n}\rangle| \le ( \frac{1}{n})^{\frac{1}{2}}$ and $ \lim_{n \to \omega} |\langle z_{n}\otimes a_{n}, x_{n}\rangle| \le  \lim_{n \to \omega} (\frac{1}{n})^{\frac{1}{2}}=0$. \vspace{\baselineskip}

Finally, for any $(z_{n})_{n}\in (M^{\omega})_{1}$ and $(w_{n})_{n}\in (N^{\omega})_{1}$, since $x_{n}\in M\bar{\otimes}A$, we have $\langle (z_{n})_{n}\otimes (w_{n})_{n}, x \rangle= \lim_{n \to \omega} \langle z_{n}\otimes w_{n}, x_{n}\rangle=\lim_{n \to \omega} \langle z_{n}\otimes E_{A}(w_{n}), x_{n} \rangle=0$, where the last equality comes from the previous paragraph. This shows that $x$ is orthogonal to $M^{\omega}\bar{\otimes}N^{\omega}$, while it is a non zero element in $M'\cap (M\bar{\otimes}A)^{\omega}$. This is enough to conclude the lemma.\begin{flushright}$\Box$
\end{flushright}

\begin{rem}
The above proof can be slightly modified to show the following statement: whenever $M,N$ are tracial von Neumann algebras with $M'\cap M^{\omega}$ and $N'\cap N^{\omega}$ diffuse, $(M\bar{\otimes}N)'\cap (M\bar{\otimes}N)^{\omega}$ is strictly larger than $(M'\cap M^{\omega})\bar{\otimes}(N'\cap N^{\omega})$.
\end{rem}

\section{PROOF OF THE MAIN THEOREM} 
The first lemma is a particular case of Ocneanu's central freedom lemma. For the proof, see \cite[Lemma 15.25]{evans1998quantum} or \cite{ioana2019class}.

\begin{lem}
Let $(M,\tau)$ be a tracial von Neumann algebra and $R \subset M$ a hyperfinite $\rm{I\hspace{-.01em}I}_{1}$ factor.
Then we have $(R'\cap R^{\omega})'\cap M^{\omega} = R \vee (R'\cap M)^{\omega}.$ In particular, $R'\cap R^{\omega}$ is a $\rm{I\hspace{-.01em}I}_{1}$ factor.
\end{lem}

As mentioned in the introduction, the proofs of the next proposition and theorem use ideas from \cite{marrakchi2018stability} and  \cite{ioana2015spectral}.
See also \cite{haagerup1985new}.

\begin{prop}Let $M$ be a tracial von Neumann algebra with $M'\cap M^{\omega}$ abelian and $A=L^{\infty}(T,\mu)$ for some probability
measure space $(T,\mu)$. Then  $M'\cap(M\bar{\otimes}A)^{\omega}$ is abelian. \label{prop3.1}   \label{proposition}
\end{prop}

\textit{Proof.} Let $(x_{n})_{n}$ and $(y_{n})_{n}$ be elements of  $(M'\cap(M\bar{\otimes}A)^{\omega})_{1}$ with $\parallel x_{n} \parallel, \parallel y_{n} \parallel \le 1$ for all $n \in \mathbb{N}$. We will show that for each  $\varepsilon>0$, we can take $V\in \omega$ such that for each $n\in V$, we have $\parallel [ x_{n} , y_{n}]\parallel_{2} <\varepsilon$.
\vspace{\baselineskip}

Since $M'\cap M^{\omega}$ is abelian, for each $\varepsilon >0$, we can take a finite subset $F\subset M$ and a $\delta>0$ such that for any $x,y\in (M)_{1}$ which satisfies $ \sum_{z\in F} \parallel [x,z]\parallel_{2}^{2} < \delta$ and  $ \sum_{z\in F} \parallel [y,z]\parallel_{2}^{2} < \delta$ , we have $\parallel [x,y]\parallel_{2}^{2} <\varepsilon_{1}=\frac{\varepsilon}{9}$.
Since $(x_{n})_{n}$ and $(y_{n})_{n}$ commute with $M$, we can take a $V\in \omega$ such that for any $n\in V$ we have 
$ \sum_{z\in F} \parallel [x_{n},z]\parallel_{2}^{2} < \delta\varepsilon_{1}$ and $ \sum_{z\in F} \parallel [y_{n},z]\parallel_{2} ^{2}< \delta\varepsilon_{1}$.\vspace{\baselineskip}

Take any $n\in V$. Consider $x_{n}$ and $y_{n}$ as $M$ valued functions on $T$ so that the above conditions translate into 
\[\displaystyle \int_{T}\sum_{z\in F}\parallel[x_{n}(t),z]\parallel_{2}^{2} dt <\delta\varepsilon_{1} ,\int_{T}\sum_{z\in F}\parallel[y_{n}(t),z]\parallel_{2}^{2} dt <\delta\varepsilon_{1} .\]

Set $X=\{t\in T|\sum_{z\in F} \parallel[x_{n}(t),z]\parallel_{2}^{2}\ge \delta\}$ and $Y=\{t\in T|\sum_{z\in F} \parallel[y_{n}(t),z]\parallel_{2}^{2}\ge \delta\}$. From the above inequalities, we have $\mu(X)<\varepsilon_{1}$ and  $\mu(Y)<\varepsilon_{1}$, so $\mu(X\cup Y)<2\varepsilon_{1}$. By the definition of $F$ and $\delta$, for $t\in(X\cup Y)^{c}$, we have $\parallel [x_{n}(t),y_{n}(t)] \parallel_{2}^{2}<\varepsilon_{1}$. Since $\parallel [x_{n}(t),y_{n}(t)]\parallel_{2}^{2}\le 4$ for all $t\in T$, we obtain
\[ \displaystyle \parallel [x_{n},y_{n}]\parallel_{2}^{2}=\int_{T} \parallel [x_{n}(t),y_{n}(t)] \parallel_{2}^{2} dt=\int_{(X\cup Y)^{c}} +\int_{X\cup Y}<\varepsilon_{1} + 4\cdot 2\varepsilon_{1}=9\varepsilon_{1}=\varepsilon .\]
This holds for any $n\in V$, which shows that this $V\in \omega$ satisfies the properties mentioned in the first paragraph. This proves that $(x_{n})_{n}$ and $(y_{n})_{n}$ commute, which shows that  $M'\cap(M\bar{\otimes}A)^{\omega}$ is abelian.\begin{flushright}$\Box$
\end{flushright}

\begin{thm}Let $M,N_{1},N_{2}$ be tracial von Neumann algebras. Assume that $M'\cap M^{\omega}$ is abelian.
Then $M'\cap(M\bar{\otimes}N_{1})^{\omega}$ and  $M'\cap(M\bar{\otimes}N_{2})^{\omega}$ commute in  $(M\bar{\otimes}N_{1}\bar{\otimes}N_{2})^{\omega}$.   \label{theorem}
\end{thm}
\textit{Proof.} Take any $(x_{n})_{n}\in (M'\cap(M\bar{\otimes}N_{1})^{\omega})_{1}$ and $(y_{n})_{n}\in (M'\cap(M\bar{\otimes}N_{2})^{\omega})_{1}$ with $\parallel x_{n} \parallel, \parallel y_{n} \parallel \le 1$ for all $n\in \mathbb{N}$. We will show that these elements commute.\vspace{\baselineskip}

Fix an orthonormal basis $\{b_{i}\}_{i}\subset N_{1}$ (respectively $\{d_{j}\}_{j} \subset N_{2}$) of $L^{2}(N_{1})$ (respectively $L^{2}(N_{2})$) and write $x_{n}=\sum_{i} a_{i}^{n}\otimes b_{i}$ and $y_{n}=\sum_{j} c_{j}^{n}\otimes d_{j}$. We may assume that these are finite sums. Since the norm of $(x_{n})_{n}$ and $(y_{n})_{n}$ are bounded by $1$, we have 
\begin{equation}
\sum_{i} a_{i}^{n}a_{i}^{n*}, \sum_{i} a_{i}^{n*}a_{i}^{n},
\sum_{j} c_{j}^{n}c_{j}^{n*}, \sum_{j} c_{j}^{n*}c_{j}^{n}\le 1 \label{a1}.
\end{equation}

Let $T_{1}$ and $T_{2}$ be copies of $\mathbb{T}^{\mathbb{N}}$ equipped with the infinite product of the usual Lebesgue probability measure. Let $u_{i}\in \mathcal{U}(L^{\infty}(T_{1}))$ (respectively $v_{j}\in \mathcal{U}(L^{\infty}(T_{2})$) be the canonical generator of the $i$th copy of $L^{\infty}(\mathbb{T})$ in $L^{\infty}(T_{1})$ (respectively $j$th copy in $L^{\infty}(T_{2})$). We set $\Theta_{n}^{1}=\sum_{i} a_{i}^{n}\otimes u_{i} \in M\bar{\otimes}L^{\infty}(T_{1})$ and $\Theta_{n}^{2}=\sum_{j} c_{j}^{n}\otimes v_{j} \in M\bar{\otimes}L^{\infty}(T_{2})$. We regard both of them as elements of $M\bar{\otimes}L^{\infty}(T_{1})\bar{\otimes}L^{\infty}(T_{2})(=\tilde{M})$. A direct computation shows the following equations.
\begin{align*}&\bullet \parallel \Theta_{n}^{1} \parallel_{2}=\parallel x_{n}\parallel_{2}, \parallel \Theta_{n}^{2} \parallel_{2}=\parallel y_{n}\parallel_{2}.\\
&\bullet  \parallel [\Theta_{n}^{1},\Theta_{n}^{2}]\parallel_{2}=\parallel [x_{n},y_{n}]\parallel_{2}.
\end{align*}

We will show that $ \lim_{n \to \omega}\parallel [\Theta_{n}^{1},\Theta_{n}^{2}]\parallel_{2}=0$, which proves the theorem.\vspace{\baselineskip}

First, we show that $U_{m}=(\Theta_{n}^{m})_{n}\in L^{2}(M'\cap \tilde{M}^{\omega}) , m=1,2$.
For all $n\in \mathbb{N}$, using (\ref{a1}) and the fact that $u_{i}$ and $v_{j}$ are orthogonal for all $i,j$, we can compute 
\begin{align*}
\tau( |\Theta_{n}^{1}|^{4})&=\sum_{i,j,k,l}\tau(a_{i}^{n*}a_{j}^{n}a_{k}^{n*}a_{l}^{n}\otimes u_{i}^{*}u_{j}u_{k}^{*}u_{l})\\
&=\sum_{i,k}\tau( a_{i}^{n*}a_{i}^{n}a_{k}^{n*}a_{k}^{n})+\sum_{i,j} \tau(a_{i}^{n*}a_{j}^{n}a_{j}^{n*}a_{i}^{n})-\sum_{i} \tau(a_{i}^{n*}a_{i}^{n}a_{i}^{n*}a_{i}^{n})\\
& \le 2.
\end{align*}
Let $p_{\lambda}^{n}= E_{[0,\lambda)}(|\Theta_{n}^{1}|)$ for $n\in \mathbb{N}$ and $\lambda> 0$. It follows from the above computation that 
\begin{equation}
\lambda^{2}\parallel \Theta_{n}^{1}(1-p_{\lambda}^{n})\parallel_{2}^{2} =\lambda^{2}\tau(|\Theta_{n}^{1}|^{2}(1-p_{\lambda}^{n})) \le \tau(|\Theta_{n}^{1}|^{4}) \le 2.
\end{equation}
This shows that $(\Theta_{n}^{1})_{n}\in L^{2}(\tilde{M})^{\omega}$ is a $\parallel \cdot \parallel_{2} $-limit point of $(\Theta_{n}^{1}p_{\lambda}^{n})_{n}\in (\tilde{M}^{\omega})_{\lambda}$ (as $\lambda \to \infty$), hence $(\Theta_{n}^{1})_{n}\in L^{2}(\tilde{M}^{\omega})$. It is easy to show that $U_{1}$ commutes with $M$ (since $(x_{n})_{n}$ commutes with $M$), so we have $U_{1} \in  L^{2}(M'\cap \tilde{M}^{\omega})$.
The same argument holds for $U_{2}$. \vspace{\baselineskip}

Using Proposition \ref{prop3.1} to $\tilde{M}$, we have that $M'\cap \tilde{M}^{\omega}$ is abelian. Since we have seen that $U_{1},U_{2}\in L^{2}(M'\cap \tilde{M}^{\omega})$, we have $U_{1}U_{2}=U_{2}U_{1}$ as elements in $L^{1}(M'\cap \tilde{M}^{\omega})$.\vspace{\baselineskip}

Computing similarly to the above, using (\ref{a1}) and that $\sum_{j,k,l,r} a_{j}^{n} a_{k}^{n*} a_{l}^{n} a_{r}^{n*}\otimes u_{j}u_{k}^{*}u_{l}u_{r}^{*}= \Theta_{n}^{1}\Theta_{n}^{1*}\Theta_{n}^{1}\Theta_{n}^{1*} \ge 0$, for all $n\in \mathbb{N}$ and $\lambda> 0$ we have
\begin{align*}
\lambda^{2} \parallel \Theta_{n}^{2} \Theta_{n}^{1}(1-p_{\lambda}^{n}) \parallel_{2}^{2} &=
\lambda^{2} \tau ( \Theta_{n}^{2} \Theta_{n}^{1} (1-p_{\lambda}^{n}) \Theta_{n}^{1*} \Theta_{n}^{2*})\\
&\leq \tau( \Theta_{n}^{2} \Theta_{n}^{1} \Theta_{n}^{1*} \Theta_{n}^{1} \Theta_{n}^{1*} \Theta_{n}^{2*})\\
&=\sum_{i,j,k,l,r,s} \tau (c_{i}^{n} a_{j}^{n} a_{k}^{n*} a_{l}^{n} a_{r}^{n*}c_{s}^{n*} \otimes v_{i}u_{j}u_{k}^{*}u_{l}u_{r}^{*}v_{s}^{*})\\
&=\tau ((\sum_{i} c_{i}^{n*}c_{i}^{n} \otimes 1)(\sum_{j,k,l,r} a_{j}^{n} a_{k}^{n*} a_{l}^{n} a_{r}^{n*}\otimes u_{j}u_{k}^{*}u_{l}u_{r}^{*}))\\
&\le \tau (\sum_{j,k,l,r} a_{j}^{n} a_{k}^{n*} a_{l}^{n} a_{r}^{n*}\otimes u_{j}u_{k}^{*}u_{l}u_{r}^{*})\\
&\le 2.
\end{align*}

Let $q_{\lambda}^{n} =E_{[0,\lambda)}(|\Theta_{n}^{2*}|)$  for $n\in \mathbb{N}$ and $\lambda> 0$. By similar computation, we have $\lambda^{2} \parallel (1-q_{\lambda}^{n})\Theta_{n}^{2} \Theta_{n}^{1} \parallel_{2}^{2} \le 2$ for all $n\in \mathbb{N}$ and $\lambda > 0$.
Hence, approximating $(\Theta_{n}^{2} \Theta_{n}^{1})_{n}$ by $(q_{\lambda}^{n} \Theta_{n}^{2} \Theta_{n}^{1} p_{\lambda}^{n})_{n}\in (M'\cap \tilde{M}^{\omega})_{\lambda^{2}}$,  we have $(\Theta_{n}^{2} \Theta_{n}^{1})_{n} \in L^{2}(M'\cap \tilde{M}^{\omega})$. The same holds for $ (\Theta_{n}^{1} \Theta_{n}^{2})_{n}$.\vspace{\baselineskip}

Now we have that $ (\Theta_{n}^{1} \Theta_{n}^{2})_{n}-(\Theta_{n}^{2} \Theta_{n}^{1})_{n}=[U_{1},U_{2}]$ is an element of $ L^{2}(M'\cap \tilde{M}^{\omega})$ and is $0$ as an element of $L^{1}(M'\cap \tilde{M}^{\omega})$. This means that $[U_{1},U_{2}]$ is also $0$ as an element of $ L^{2}(M'\cap \tilde{M}^{\omega})$. This shows that $\lim_{n \to \omega} \parallel  [\Theta_{n}^{1},\Theta_{n}^{2}]\parallel_{2}=\parallel [U_{1},U_{2}]\parallel_{2}=0$, as we wanted.
\begin{flushright}$\Box$
\end{flushright}

\begin{cor}Let $M$ be a non-McDuff $\rm{I\hspace{-.01em}I}_{1}$ factor and $N$ any  $\rm{I\hspace{-.01em}I}_{1}$ factor. Suppose we have $M\bar{\otimes}N=K\bar{\otimes}R(=\tilde{M})$ where $K$ is a  $\rm{I\hspace{-.01em}I}_{1}$ factor and $R$ is the hyperfinite  $\rm{I\hspace{-.01em}I}_{1}$ factor. Then we have $R'\cap R^{\omega}\subset(K\vee N)^{\omega}$. \label{cor3.4}
\end{cor}
\textit{Proof.} Let $N_{1},N_{2}$ be copies of $N$. Let $M\bar{\otimes}N_{1}=K_{1}\bar{\otimes}R_{1}$ and $M\bar{\otimes}N_{2}=K_{2}\bar{\otimes}R_{2}$ be the corresponding tensor product decomposition where $K\cong K_{1}\cong K_{2}$ and $R\cong R_{1}\cong R_{2}$. We regard both of them as  subalgebras of $M\bar{\otimes}N_{1}\bar{\otimes}N_{2}$.\vspace{\baselineskip}

From the previous theorem, we have that $ M'\cap (M\bar{\otimes}N_{1})^{\omega}$ and $M'\cap (M\bar{\otimes}N_{2})^{\omega}$ commute. Since  we have $R_{1}'\cap R_{1}^{\omega} \subset M'\cap (M\bar{\otimes}N_{1})^{\omega}$ and $ R_{2}'\cap R_{2}^{\omega} \subset M'\cap (M\bar{\otimes}N_{2})^{\omega}$, we obtain $R_{1}'\cap R_{1}^{\omega} \subset (R_{2}'\cap R_{2}^{\omega})'\cap (M\bar{\otimes}N_{1}\bar{\otimes}N_{2})^{\omega}$. Using Lemma 3.1 to $R_{2} \subset M\bar{\otimes}N_{1}\bar{\otimes}N_{2}$, we obtain $R_{1}'\cap R_{1}^{\omega} \subset R_{2} \vee (K_{2}\bar{\otimes}N_{1})^{\omega}$. Since $R_{1}'\cap R_{1}^{\omega}$ commutes with the whole $M\bar{\otimes}N_{1}\bar{\otimes}N_{2}$, we moreover have $R_{1}'\cap R_{1}^{\omega} \subset (K_{2}\bar{\otimes}N_{1})^{\omega}\cap  (M\bar{\otimes}N_{1})^{\omega}$.\vspace{\baselineskip} 

Take any $(\sum_{i} x_{i}^{n}\otimes y_{i}^{n})_{n}\in R_{1}'\cap R_{1}^{\omega} \subset (K_{2}\bar{\otimes}N_{1})^{\omega}\cap  (M\bar{\otimes}N_{1})^{\omega}$, where $x_{i}^{n}\in K_{2}$ and $y_{i}^{n}\in N_{1}$. We may assume that the sum is a finite sum.
Since it is in $(M\bar{\otimes}N_{1})^{\omega}$, we have $(\sum_{i} x_{i}^{n}\otimes y_{i}^{n})_{n}=(\sum_{i} E_{M}(x_{i}^{n})\otimes y_{i}^{n})_{n}$.
Denote by $\sigma \in \text{Aut}(M\bar{\otimes}N_{1}\bar{\otimes}N_{2})$ the automorphism which flips $N_{1}$ and $N_{2}$. We have $\sigma (K_{2}) = K_{1}$ and $E_{M} \circ \sigma =E_{M}$. Thus, $(\sum_{i} x_{i}^{n}\otimes y_{i}^{n})_{n}= (\sum_{i} E_{M}(\sigma(x_{i}^{n}))\otimes y_{i}^{n})_{n}$ with $\sigma(x_{i}^{n})\in K_{1}$. Since $M=N_{1}'\cap (M\bar{\otimes}N_{1})$, we have $E_{M}(\sigma(x_{i}^{n}))\in \overline{\text{co}}^{\text{w.o}.}(\{u\sigma (x_{i}^{n})u^{*}|u\in \mathcal{U}(N_{1})\})$ (the w.o. closed convex hull of $\{u\sigma (x_{i}^{n})u^{*}|u\in \mathcal{U}(N_{1})\}$ ) and thus $\sum_{i} E_{M}(\sigma(x_{i}^{n}))\otimes y_{i}^{n} \in K_{1}\vee N_{1}$. This shows 
 $(\sum_{i} x_{i}^{n}\otimes y_{i}^{n})_{n}\in (K_{1}\vee N_{1})^{\omega}$, hence  $R_{1}'\cap R_{1}^{\omega} \subset  (K_{1}\vee N_{1})^{\omega}$. This proves the corollary.
\begin{flushright}$\Box$
\end{flushright}
Note that by \cite{chifan2018remark} (or \cite{ioana2019class}),  $R'\cap R^{\omega}\subset(K\vee N)^{\omega}$ implies $R\prec_{M\bar{\otimes}N}K\vee N$, that is, $R$ embeds into $K\vee N$ in $M\bar{\otimes}N$, in the sense of Popa's intertwining technique \cite{popa2003strrgd}. Since $K \subset K\vee N$ and that $E_{K\vee N}$ is a $K$-bimodule map, we moreover have $M\bar{\otimes}N=K\bar{\otimes}R\prec_{M\bar{\otimes}N}K\vee N$. This gives information on the position of $K$ in $M\bar{\otimes}N$ whenever $M\bar{\otimes}N=K\bar{\otimes}R$ is a tensor product decomposition of $M\bar{\otimes}N$ with $R$ hyperfinite. In the context of McDuff decomposition, this gives the following statement.

\begin{cor} Let $M$ be a non-McDuff $\rm{I\hspace{-.01em}I}_{1}$ factor and $R_{1}\cong R$ . Set $\tilde{M}=M\bar{\otimes}R_{1}$. Suppose we have another tensor product decomposition $\tilde{M}=K\bar{\otimes}R_{2}$ with $K$ a $\rm{I\hspace{-.01em}I}_{1}$ factor and $R_{2}\cong R$. Then we have $\tilde{M}=M\bar{\otimes}R_{1}\prec_{\tilde{M}} K\vee R_{1}$. \label{cor3.5}
\end{cor}

Note that by \cite[Theorem B]{hoff2016neumann}, we can not prove $M\prec_{\tilde{M}}K$ in this setting which allows $M$ to have property Gamma. However, the above corollary shows that if we consider $M\bar{\otimes}R_{1}$ and $K\vee R_{1}$ instead of $M$ and $K$, we always have $\tilde{M}=M\bar{\otimes}R_{1}\prec_{\tilde{M}} K\vee R_{1}$.

\section{ACKNOWLEDGEMENTS} 
This paper is written for master's thesis of the author. First of all, the author would like to thank his supervisor, Professor Yasuyuki Kawahigashi for his helpful comments and continuing support. The author is  grateful to Yusuke Isono and Narutaka Ozawa for their comments which have greatly improved this paper. The author is also grateful to the anonymous referee for several suggestions which led to Corollary \ref{cor3.5} which was not included in the original version of this paper.

\bibliographystyle{amsalpha}
\bibliography{myref(1)}

\small{\textsc{Graduate School of Mathematics, University of Tokyo, Komaba, Tokyo 153-8914, Japan}
\textit{E-mail adress:} y.hashiba00@gmail.com}

\end{document}